\documentclass[11pt]{amsart}%
\usepackage{amssymb, amsmath}
\usepackage{xcolor}
\usepackage{graphicx}
\usepackage{hyperref}
\usepackage{geometry}
\usepackage{amsmath}
\usepackage{amsfonts}
\usepackage{amssymb}%
\definecolor{cmk2020}{HTML}{BE1E2D}
\hypersetup{
    colorlinks=true,
    citecolor=cmk2020,
    urlcolor=cmk2020
}
\newtheorem{theorem}{Theorem}[section]

\theoremstyle{definition}

\newtheorem{example}[theorem]{Example}

\newtheorem{corollary}[theorem]{Corollary}
\theoremstyle{remark}

\numberwithin{equation}{section}
\begin{document}
\title{Central Open Sets Tilings}
\author[L.F. Barnsley]{Louisa F. Barnsley}
\author[M.F. Barnsley]{Michael F. Barnsley}

\begin{abstract}
We introduce a method for constructing collections of subsets of
$\mathbb{R}^{n}$, using an iterated function system, a set $T,$ and a cost
function. We refer to these collections as tilings. The special case where $T$
is the central open set of an iterated function system that obeys the open set
condition is emphasized. The notion of the central open set associated with an
iterated function system of similitudes, introduced in 2005 by Bandt, Hung,
and Rao, is reviewed. A practical method for calculating pictures of central
open sets is described. Some general properties and examples of the tilings
are presented.

\vspace{0.1cm}

\end{abstract}
\maketitle

\baselineskip12pt



\section{Introduction}

\begin{figure}[ptb]
\centering
\includegraphics[
height=1.9867in,
width=4.846in
]{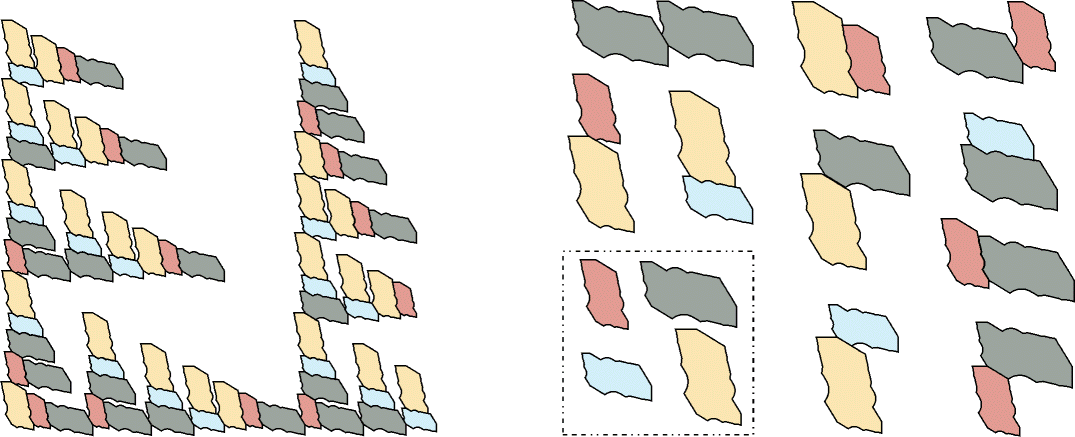}\caption{See Example \ref{combplusex}. The left panel
illustrates part of an unbounded central open set tiling. The right panel
shows some of the ways in which the four prototiles (see inset box) can meet.
See also Figures \ref{combplus4} and \ref{combtilefill4}$.$}%
\label{exactsbys}%
\end{figure}

The goal of this paper is to describe a simple method for producing a wide
range of tilings (in a generalized sense). The method uses an iterated
function system (IFS) acting on $\mathbb{R}^{n}$, together with a set $T$, and
a cost function $c$. Figure \ref{exactsbys} illustrates part of such a tiling
where $T$ is the central open set of an IFS. The possible structures of the
new tilings are diverse, yet all are handled with the same underlying
mathematical device. The formalism yields examples in analysis, geometry, and
dynamics; it leads to extensive and rich families of variants of the
self-similar tilings introduced in \cite{polygon, barnsleyvince}. Our methods
are based on addresses associated with IFSs and mappings from these addresses
into tilings and tiling spaces.

Our main examples use the central open set of the attractor of an IFS to
provide the shapes of the tiles. In this case the tilings have properties that
suggest they may be used to model patterns that arise naturally. For example,
such tiles have self-similar features and may touch, and there may also be
gaps that repeat at different scales. Type \textquotedblleft
mudcracks\textquotedblright\ into a search engine to see illustrations of
seemingly related natural patterns.

In this paper we use the following lexicon. We use the word \textit{tiling}%
\ to mean a collection of closed subsets of $\mathbb{R}^{n}$, and a
\textit{tile}\ is a member of the collection. This is more general than the
standard definitions, see Grunbaum and Sheppard \cite{grun}. We say that two
\textit{tiles} \textit{meet }if their intersection is non-empty, and we say
that they \textit{touch} if their intersection is non-empty and contains no
interior points. The \textit{support }of a tiling is the union of its tiles.
Two \textit{tilings} \textit{meet} if the union of the tiles in their
intersection is a tiling whose support is the intersection of the supports of
the two tilings. A set of \textit{prototiles} for a tiling is a set of tiles
such that each member of the tiling is related to a prototile by an isometry.
The set of isometries may be restricted to translations. We say that two tiles
have the same \textit{shape} if they are related by an isometry. If two tiles
are related by an isometry we say that they are \textit{copies }of one
another. We say that a tiling is \textit{commensurate }if the sizes of all its
tiles belong to a geometrical progression; otherwise the tiling is
\textit{incommensurate;} see also \cite{smilansky}\textit{. }A \textit{patch}
of a tiling is the set of its tiles that have non-empty intersection with a
finite set, typically a disk or rectangle.

In this paper we are concerned with the situation where all tiles have the
same shape. But the theory is readily generalized to multiple shapes, by using
graph-directed IFS, as in \cite{barnvince2020}.

\section{Iterated function systems and central open sets}

Here we review the notion of an IFS of similitudes and of a central open set
introduced in 2005 by Bandt, Hung, and Rao, \cite{bandtneighbor2}.

Roughly following \cite{bandtneighbor2}, let $F=\{\mathbb{R}^{n};f_{1}%
,f_{2}\dots f_{m}\}$ denote a collection of contractive similitudes, that is
$f_{i}:\mathbb{R}^{n}\rightarrow\mathbb{R}^{n}$ with%
\[
|f_{i}(x)-f_{i}(y)|=\lambda_{i}|x-y|\text{ for all }x,y\in\mathbb{R}^{n}%
\]
where the $\lambda_{i}\in(0,1)$ are the contraction factors and $\left\vert
\cdot\right\vert $ denotes the Euclidean norm. We refer to $F$ as an iterated
function system. By slight abuse of notation we use the same symbol $F$ to
denote the mapping from sets to sets $F:2^{\mathbb{R}^{n}}\rightarrow$
$2^{\mathbb{R}^{n}}$ defined by $F(S)=\{f_{i}(x):x\in S,i=1,\dots m\}.$ It is
well-known that there exists a unique non-empty compact set $A\subset
\mathbb{R}^{n}$ such that
\begin{align*}
A  &  =F(A)=\bigcup\limits_{i=1}^{m}f_{i}(A)\\
&  \text{ }%
\end{align*}
where $f_{i}(A)=\{f_{i}(x)|x\in A\}$, \cite{hutchinson}$.$

The set $A$ is called the \textit{attractor} of the iterated function system
$F,$ because
\[
\lim_{k\rightarrow\infty}F^{k}(\{x\})=A\text{ for all }x\in\mathbb{R}^{n},
\]
where convergence is with respect to the Hausdorff metric on $\mathbb{R}^{n}$,
and $F^{k}$ is the function $F$ composed with itself $k$ times. We say that
the \textit{basin} of $A$ is $\mathbb{R}^{n}$. The \textit{fast basin}
\cite{barnlesniak, barnfast} of $A$ is a subset of the basin defined by
\[
B=\{x\in\mathbb{R}^{n}|F^{k}(\{x\})\cap A\neq\emptyset\text{, some }%
k\in\mathbb{N}\},
\]
where $\mathbb{N}$ is the set of positive integers. The fast basin is the set
of points such that some finite orbit meets $A.$ Fast basins are related to
but distinct from the fractal blow-ups introduced by Strichartz
\cite{strichartz} and the macro-fractals\textit{ }introduced by Banakh and
Novosad \cite{banakh}. We will use $B\backslash A$ in calculations in Section
4$.$

Roughly quoting \cite{bandtneighbor2}, the attractor $A$ is the union of
smaller copies of itself, $A_{i}=f_{i}(A),$ where each $A_{i}$ consists of
smaller copies $A_{ij}=f_{i}(f_{j}(A)),$ and so on. For any positive integer
$k,$ we can consider the set $\Sigma^{k}$ of words $\mathbf{i=}i_{1}\dots
i_{k}$ from the alphabet $\Sigma=\{1,2,\dots m\}$. Writing $f_{\mathbf{i}%
}=f_{i_{i}}f_{i_{2}}\dots f_{i_{k}}$ and $A_{\mathbf{i}}=f_{\mathbf{i}}(A)$ we
have%
\[
A=\bigcup\{f_{\mathbf{i}}(A)|\mathbf{i\in}\Sigma^{k}\}\text{.}%
\]
When $k$ tends to infinity, this induces a continuous map that we call the
\textit{address }map, $\pi:\Sigma^{\infty}\rightarrow A,$ from the set
$\Sigma^{\infty}$ of infinite sequences $i_{1}i_{2}\dots$ onto the attractor.

The IFS $F$ is said to satisfy the open set condition (OSC) if there is a
nonempty open set $O\subset\mathbb{R}^{n}$ such that%
\[
F(O)\subset O\text{ and }f_{i}(O)\cap f_{j}(O)=\emptyset\text{ for }i\neq j.
\]
Such a set $O$ is called a \textit{feasible open set} of $F$. The OSC plays an
important role in fractal geometry. For example, if $F$ obeys the OSC, then
\cite{moran} the Hausdorff dimension of $A$ is the unique positive solution
$D$ of%
\[
\sum\limits_{i=1}^{m}\lambda_{i}^{D}=1.
\]
See \cite{bandtneighbor2} for a succinct account of the history and
mathematical significance of the OSC. Here we are interested in a particular
feasible open set, the \textit{central open set} of $F$, and its relationship
to fractal tilings \cite{barnsleyvince}.

The second requirement of the OSC may be written%
\[
O\cap f_{\mathbf{i}}^{-1}f_{\mathbf{j}}(O)=\emptyset
\]
for $i_{1}\neq j_{1}$. The maps in%
\[
\mathcal{N=\{}f_{\mathbf{i}}^{-1}f_{\mathbf{j}}|\mathbf{i,j}\in\Sigma^{\ast
},i_{1}\neq j_{1}\}\text{ where }\Sigma^{\ast}=\bigcup\limits_{k=1}^{\infty
}\Sigma^{k}%
\]
are called \textit{neighbor maps}, \cite{bandtneighbor1, bandtneighbor2}.
Neighbor maps may be used to provide an algebraic formulation of the OSC:
there is a constant\textit{ }$\kappa>0$ such that\textit{ }$\left\Vert
h-id\right\Vert >\kappa$ for all neighbor maps $h.$ Neighbor maps are related
to the fast basin by
\[
B\backslash A=H\backslash A\text{ where }H=\bigcup\{h(A)|h\in\mathcal{N\}}.
\]
Any feasible open set $O$ must have empty intersection with $H$.

The \textit{central open set} $C$ for $F$ is defined to be
\[
C=\{x\in\mathbb{R}^{n}|d(x,A)<d(x,H)\}=\{x\in\mathbb{R}^{n}%
|d(x,A)<d(x,B\backslash A)\}
\]
where $d(x,Y)=\inf\{\left\vert x-y|\text{ }\right\vert y\in Y\}$. Bandt et al.
\cite{bandtneighbor2} prove the following theorem and its elegant corollary.

\begin{theorem}
If the OSC holds, then the central open set $C$ is a feasible open set. If the
OSC does not hold then $C$ is empty.
\end{theorem}

\begin{corollary}
The OSC holds if and only if $A$ is not contained in $\overline{H}$.
\end{corollary}

It is an unanswered question as to whether or not it is true that the OSC
holds if and only if $A\neq\overline{A\cap H}$. The latter was claimed by M.
Moran \cite{moran}, but his proof contains a gap \cite{bandtneighbor2}$.$

Questions relating to the existence of, and the structure of, a feasible open
set for a given IFS are very subtle, see for example \cite{tetenov}%
.\begin{figure}[ptb]
\centering\includegraphics[
height=1.5849in,
width=5.1067in
]{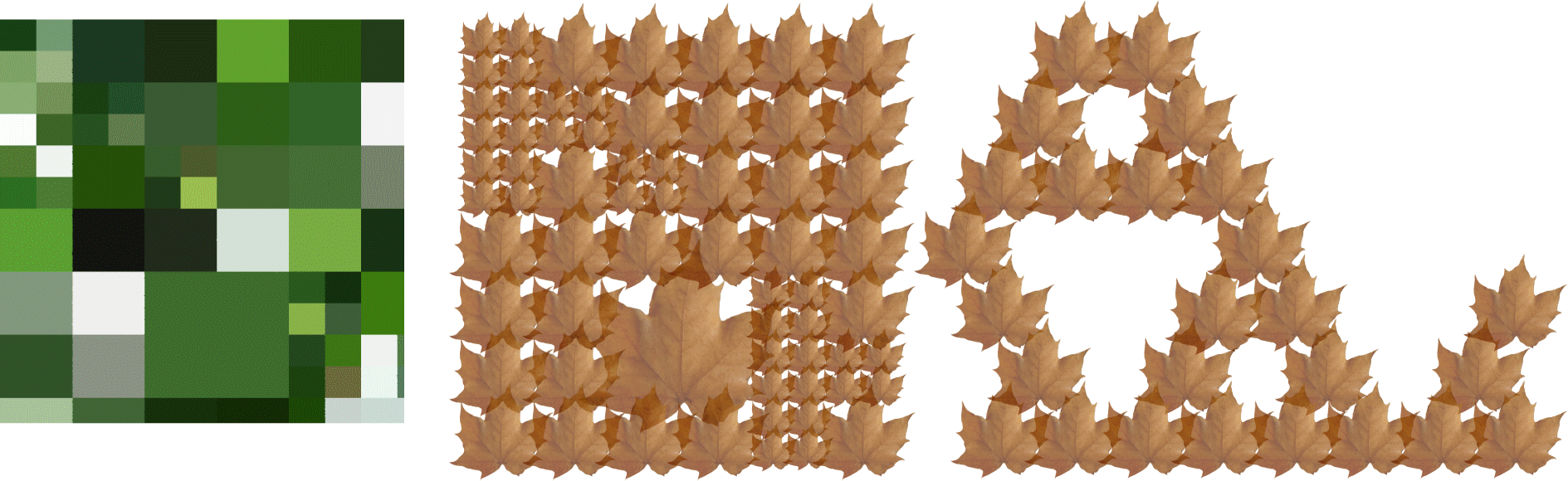}\caption{Examples of patches of three unbounded tilings. The
leftmost two images are related to an IFS whose attractor is a square. The
rightmost image is related to an IFS whose attractor is a Sierpinski triangle.
Two of the panels use $T=L,$ a decorated leaf-shaped set,  see Example
\ref{leaftilesex}.}%
\label{leaftiles4}%
\end{figure}

\section{Diverse tilings derived from an IFS and a cost function}

In this Section we describe a general construction of tilings using an IFS
$F$, a set $T\subset\mathbb{R}^{n}$, and a cost function $c$ defined below.
The resulting tilings may have overlapping tiles with non-empty interiors.
Such tilings might be used to model fallen leaves carpeting a forest floor,
duckweed on the surface of a pond, cracks in dried mud, or to design patterns
for wallpaper.

Let $F$ be a IFS consisting of at least two distinct similitudes, and let
$T\ $be a closed subset of $\mathbb{R}^{n}$. For convenience we suppose that
$A\cap T\neq\emptyset,$ but this is not necessary. We use sets of similitudes,
which are subsets of $\mathcal{N}$ determined by the cost function and are
applied to $T$, to form collections of scaled and translated and possibly
flipped (i.e. turned upside down, in the two-dimensional case) copies of $T,$
with possible overlaps.

Let $\mathbb{H(R}^{n})$ be the closed bounded subsets of $\mathbb{R}^{n}$
equipped with the spherical Hausdorff metric$,$ see \cite{barnvince2019,
barnvince2020}. This metric $d$ is defined as follows. Let $P(x)$ be the
stereographic projection of $x\in\mathbb{H(R}^{n})$ onto the $\left(
n+1\right)  $ dimensional sphere tangent to $\mathbb{R}^{n}$ at the origin.
Then $d(x,y)$ is the Hausdorff distance between $P(x)$ and $P(y)$ using the
round metric $d_{R}$ on the sphere. Let $d_{H}(X,Y)$ be the Hausdorff distance
between sets of subsets $X$ and $Y$ of $\mathbb{R}^{n}$ calculated using
$d_{R}.$ Let $\mathbb{H(H(R}^{n}))$ be the collections of subsets of
$\mathbb{H(R}^{n})$ that are closed with respect to $d_{H}$.

For each $i\in\left\{  1,\dots M\right\}  $ write $\lambda_{i}=s^{a_{i}}$
where $s=\max\{\lambda_{i}|i=1,...m\}$ and assign a cost $c_{i}>0$ to the map
$f_{i}.$ For example we may choose $c_{i}=a_{i}$. For $\mathbf{i}%
\mathbf{=}i_{1}i_{2}\dots\in\Sigma^{\infty},$ write $\mathbf{i}|k=i_{1}\dots
i_{k}\in\Sigma^{\ast},$ and $\mathbf{i}|0=\emptyset$. Define a \textit{cost
function} $c:\{\emptyset\}\cup\Sigma^{\ast}\rightarrow(0,\infty)$ by
\[
c(\mathbf{i}|k)=c_{i_{1}}+c_{i_{2}}+\dots+c_{i_{k}}\text{, }c(\emptyset)=0.
\]
Define a mapping $\Pi_{T}:\{\emptyset\}\cup\Sigma^{\ast}\cup\Sigma^{\infty
}\rightarrow\mathbb{H(H(R}^{n}))$ by%
\begin{align*}
\Pi_{T}(\mathbf{i}|k) &  =f_{-(\mathbf{i}|k)}(\{f_{(\mathbf{j}|l)}%
(T)|\mathbf{j}\in\Sigma^{\infty},l\in\mathbb{N},\text{ }c(\mathbf{j}|l-1)\leq
c(\mathbf{i}|k)<c(\mathbf{j}|l)\}),\\
\text{ }f_{-(\mathbf{i}|k)} &  :=f_{i_{1}}^{-1}\dots f_{i_{k}}^{-1}%
,f_{(\mathbf{j}|l)}=f_{j_{1}}...f_{j_{l}},\\
\Pi_{T}(\mathbf{i}) &  =\bigcup\limits_{k=1}^{\infty}\Pi_{T}(\mathbf{i}%
|k)\text{ for all }\mathbf{i}\in\Sigma^{\infty},\Pi_{T}(\emptyset
)=\{f_{1}(T),f_{2}(T),...f_{m}(T)\}.
\end{align*}
$\Pi_{T}$ is well-defined because $\left\{  \Pi_{T}(\mathbf{i}%
|k)|k=1,2,...\right\}  $ is a nested increasing sequence of collections of
sets:
\begin{equation}
\Pi_{T}(\mathbf{i}|0)\subset\Pi_{T}(\mathbf{i}|1)\subset\Pi_{T}(\mathbf{i}%
|2)\dots\label{nestedeq}%
\end{equation}
for all $\mathbf{i}\in\Sigma^{\infty}$. This is true because%
\begin{align*}
\Pi_{T}(\mathbf{i}|k+1) &  =f_{-(\mathbf{i}|k+1)}(\{f_{(\mathbf{j}%
|l)}(T)|\mathbf{j}\in\Sigma^{\infty},l\in\mathbb{N},c(\mathbf{j}|l-1)\leq
c(\mathbf{i}|k+1)<c(\mathbf{j}|l)\})\\
&  \supset f_{-(\mathbf{i}|k+1)}(\{f_{(\mathbf{j}|l)}(T)|\mathbf{j}\in
\Sigma^{\infty},l\in\mathbb{N},c(\mathbf{j}|l-1)\leq c(\mathbf{i}%
|k+1)<c(\mathbf{j}|l),j_{1}=i_{k+1}\})\\
&  =f_{i_{1}}^{-1}\dots f_{i_{k+1}}^{-1}(\{f_{i_{k+1}}f_{j_{2}}...f_{j_{l}%
}(T)|\mathbf{j}\in\Sigma^{\infty},l\in\mathbb{N},c(\mathbf{j}|l-1)\leq
c(\mathbf{i}|k)<c(\mathbf{j}|l)\}\\
&  =\Pi_{T}(\mathbf{i}|k).
\end{align*}
Equation (\ref{nestedeq}) is the key mathematical device in this paper. In
general $\Pi_{T}(\mathbf{i})$ is a collection of subsets of $\mathbb{R}^{n}$,
that we call tiles. These tiles are translated, scaled, maybe flipped and/or
rotated, copies of $T$. They may be overlapping and the support of the tiling
$\Pi_{T}(\mathbf{i)}$, namely $\bigcup\{t\in\Pi_{T}(\mathbf{i})\}\subset
\mathbb{R}^{n}$ may be complicated.

If $F$ obeys the OSC and $T\subset\overline{C},$ then distinct sets of the
form $f_{-(\mathbf{i}|k)}f_{(\mathbf{j}|l)}(T)$ in the tiling $\Pi
_{T}(\mathbf{i})$ are non-overlapping; that is, the interiors of the
intersections of distinct tiles are empty.

Denote the range of $\Pi_{T}:\Sigma^{\infty}\rightarrow\mathbb{H(H(R}^{n}%
))\ $by $\mathbb{T}_{T}=\{\Pi_{T}(\mathbf{i)}|\mathbf{i}\in\Sigma^{\infty
}\mathbf{\}}$. As a consequence of properties of, and structures associated
with, the shift map $\sigma:\Sigma^{\infty}\rightarrow\Sigma^{\infty},$ much
can be said, along the lines of \cite{barnvince2019, barnvince2020}, about
continuity properties of $\Pi_{T}:\Sigma^{\infty}\rightarrow\mathbb{T}_{T}$
with respect to the metric $d_{H}$ defined above, dynamics, invariant
measures, and ergodic properties associated with mappings that take
$\mathbb{T}_{T}$ into itself, such as certain inflation and deflation operations.

\begin{example}
Let $F=\mathbb{\{R}^{1};f_{1},f_{2}\},f_{1}(x)=\frac{x}{2},f_{2}(x)=\frac
{x+1}{2}.$ Then $A=[0,1]$. Choosing $T=[-\frac{1}{3},\frac{4}{3}]$, we find
$\Pi_{\lbrack-\frac{1}{3},\frac{4}{3}]}(\overline{1})=$ $\{[-\frac{1}{6}%
+\frac{n-1}{2},\frac{1}{6}+\frac{n}{2}]|$ $n=1,2,\dots\}$. That is,
$\Pi_{\lbrack-\frac{1}{3},\frac{4}{3}]}(\overline{1})$ is a collection of
overlapping closed intervals whose union is $[-\frac{1}{6},\infty)$.
\end{example}

\begin{example}
\label{leaftilesex} (i) The leftmost panel in Figure \ref{leaftiles4}
illustrates part of a tiling $\Pi_{A}(\overline{1})$ where the IFS of four
similitudes each with scaling factor $0.5$, attractor $A=[0,1]\times
\lbrack0,1]$, and the cost function is defined by $c_{1}=1,$ $c_{2}=1.3,$
$c_{3}=1.5,$ $c_{4}=$ $2.$ (ii) The middle panel illustrates the same part,
but of $\Pi_{L}(\overline{1})$ where $L$ is the support of a leaf picture. In
this case the tiles have been decorated by a picture of a leaf. (iii) The
rightmost panel is related to an IFS of three maps, whose attractor is a
Sierpinski triangle, and the cost function specified by $c_{1}=c_{2}=c_{3}=1,$
and the same set $L$.
\end{example}

For a two-dimensional affine transformation $f:\mathbb{R}^{2}\rightarrow
\mathbb{R}^{2}$ we write
\[
\text{ }f=%
\begin{bmatrix}
a & b & e\\
c & d & g
\end{bmatrix}
\text{ for }f(x,y)=(ax+by+e,cx+dy+g)\text{ }%
\]
where $a,b,c,d,e,g\in\mathbb{R}$.

\begin{figure}[ptb]
\centering\includegraphics[
height=2.0913in,
width=4.3918in
]{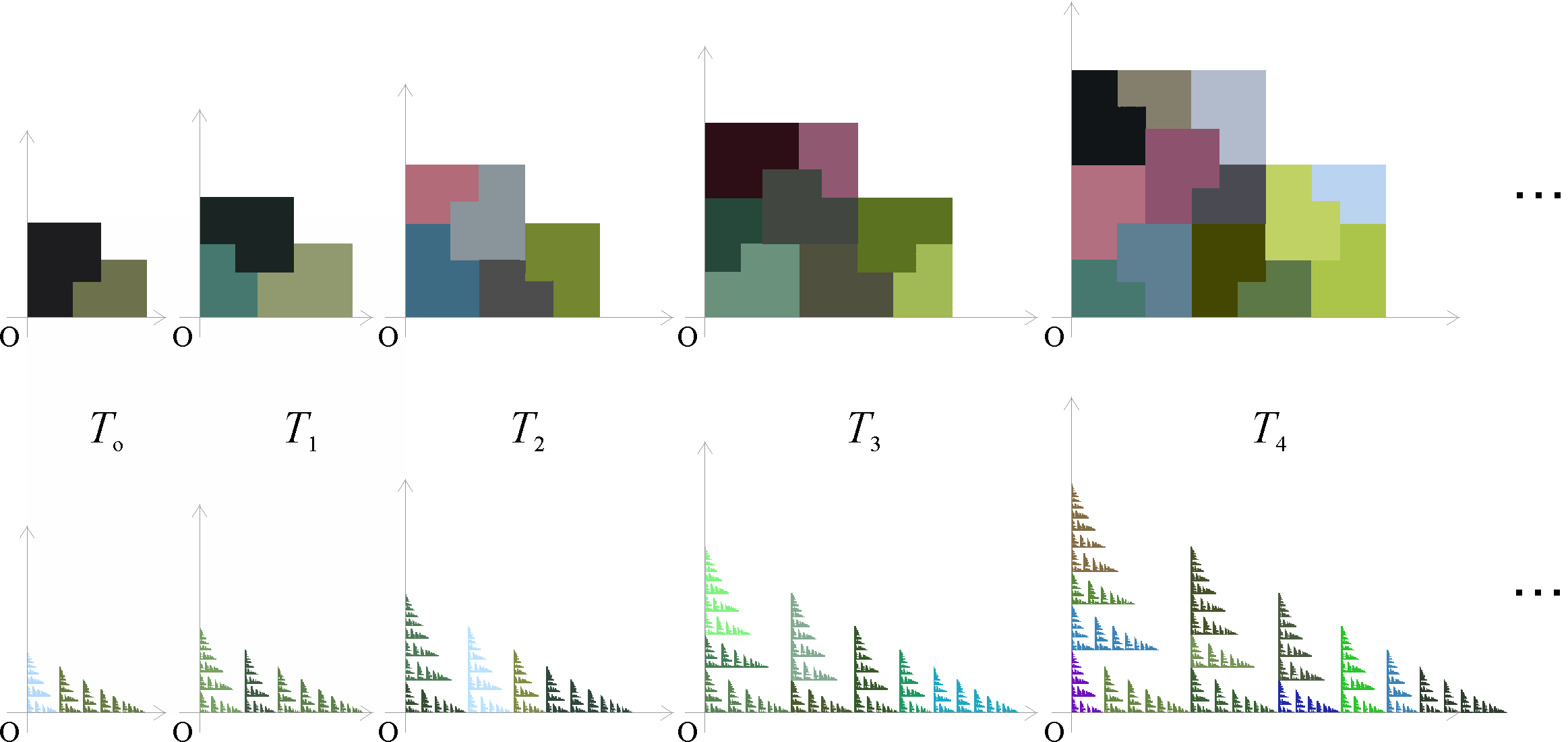}\caption{Canonical tilings $\left\{
T_{k}\right\}  $ for two different iterated function systems $\{\mathbb{R}%
^{2},f_{1},f_{2}\}$. The tiling $T_{k+1}$ is derived from $T_{k}$ by replacing
each isometric copy of $A$ in $s^{-1}T_{k}$ by a copy of $T_{0}=\{f_{1}%
(A),f_{2}(A)\}.$ See Example \ref{fig2-introex}.}%
\label{fig2-intro}%
\end{figure}

\begin{example}
\label{fig2-introex} In the special case where $c_{i}=a_{i}$ $\in\mathbb{N}$
for all $i=1,\dots m$, we call the sets
\[
T_{k}=s^{-k}\{f_{(\mathbf{i}|l)}(A)|c(\mathbf{i}|l-1)\leq k<c(\mathbf{i}%
|l)\},T_{0}=F(A),
\]
\textit{canonical tilings}. They play a natural role in fractal tilings
\cite{barnvince2020}, and in connecting them to algebraic geometry
\cite{anderson}. Two sequences of canonical tilings are illustrated in Figure
\ref{fig2-intro}. The IFSs may be deduced from $T_{0}=\{f_{1}(A),f_{2}(A)\}$.
The IFS for the top sequence is $F=\left\{  \mathbb{R}^{2};f_{1}%
,f_{2}\right\}  $ where%
\[
f_{1}=%
\begin{bmatrix}
0 & s & 0\\
-s & 0 & s
\end{bmatrix}
,f_{2}=%
\begin{bmatrix}
-s^{2} & 0 & 1\\
0 & s^{2} & 0
\end{bmatrix}
\]
where $s+s^{2}=1,$ $s>0,$ and for the lower sequence
\[
f_{1}=%
\begin{bmatrix}
s & 0 & 1\\
0 & s & 0
\end{bmatrix}
,f_{2}=%
\begin{bmatrix}
0 & s^{2} & 0\\
s^{2} & 0 & 0
\end{bmatrix}
\]
with $s+s^{4}=1$, $s>0$. In these two cases, and others like them, when the
cost function is defined by $c_{i}=i,$ there is a simple relationship between
the canonical tilings and the tilings $\Pi_{A}(\mathbf{i}),$ namely
\[
\Pi_{A}(\mathbf{i|}k)=f_{-(\mathbf{i|}k)}s^{c(\mathbf{i}|k)}T_{c(\mathbf{i}%
|k)}\text{, where }c(\mathbf{i|}k\mathbf{)=}i_{1}+i_{2}+\dots i_{k},
\]
for all $\mathbf{i}$ and $k$. That is, $\Pi_{A}(\mathbf{i|}k)$ is isometric to
$T_{c(\mathbf{i}|k)},$ see \cite{barnvince2020}.
\end{example}

\section{Central open set tilings}

In this Section we assume that $F$ obeys the OSC, and consider the two special
cases $T=A$ and $T=$ $\overline{C}$ in the mapping $\Pi_{T}:\{\emptyset
\}\cup\Sigma^{\ast}\cup\Sigma^{\infty}\rightarrow\mathbb{H(H(R}^{n})).$ For
consistency with \cite{barnvince2020} we define
\[
\Pi=\Pi_{A}\text{ and }\Xi=\Pi_{\overline{C}}\text{ .}%
\]
These tilings are of particular interest to us. Here's why. If $A$ has
nonempty interior, so that $\overline{C}=A$, these tilings may be examples of
conventional self-similar tilings as defined by \cite{grun}, or tilings with
fractal boundaries \cite{kenyon, lagarias}. But they are more general because
they may be tilings with infinitely many incommensurate tile sizes. The case
$T=$ $\overline{C}$ is special because it seems to be an extreme case: we
conjecture that if $T$ is chosen to be a closed set that contains an open set
that contains $\overline{C},$ then $\Pi_{T}(\mathbf{i})$ contains overlapping tiles.

The mapping $\Pi=\Pi_{A}$, and the fractal tilings it generates when
$c_{i}=a_{i}$, were introduced and studied in \cite{tilings, polygon,
barnsleyvince}. We refer to $\Pi(\mathbf{i})$ as a fractal tiling, and we
refer to sets of the form $f_{(\mathbf{i|}k\mathbf{)}}^{-1}f_{(\mathbf{j|}%
l\mathbf{)}}(A)$ as fractal tiles. Note that fractal tiles may have empty
interiors. They have non-empty interiors when $A$ has non-empty interior. The
relationship of the address map $\pi:\Sigma^{\infty}\rightarrow A$ to the
contractive IFS $F$ has analogies with the relationship of $\Pi$ to the
expansive IFS $F^{-1}=\{\mathbb{R}^{2};f_{1}^{-1},\dots f_{m}^{-1}\}$, see
also \cite{banakh, strichartz}.

We call $\Xi(\mathbf{i})=\Pi_{\overline{C}}(\mathbf{i}$) a \textit{central
open set tiling}. We refer to sets of the form $f_{(\mathbf{i|}k\mathbf{)}%
}^{-1}f_{(\mathbf{j|}l\mathbf{)}}(\overline{C})$ as central open set
tiles\textit{. }The interiors of the tiles in $\Xi(\mathbf{i})$ are non-empty
and disjoint for any fixed $\mathbf{i}$, but the tiles may touch.

\begin{example}
Let $F=\mathbb{\{R}^{1};f_{1},f_{2}\},f_{1}(x)=\frac{x}{2},f_{2}(x)=\frac
{x+1}{2}.$ Then $A=[0,1]$, $C=(0,1)$ and we find $\Pi(\overline{1}%
)=\Xi(\overline{1})=$ $\{[\frac{n-1}{2},\frac{n}{2}]\subset\mathbb{R}|$
$n\in\mathbb{N}\}$, and if the tail of $\mathbf{i\in}\Sigma$ is neither
$\overline{1}=11...$ nor $\overline{2}=22...,$ then $\Pi(\mathbf{i}%
)=\Xi(\mathbf{i})=\{[\frac{n-1}{2},\frac{n}{2}]|$ $n\in\mathbb{Z}\}$. If the
tail of $\mathbf{i\in}\Sigma$ is either $\overline{1}=11...$ or $\overline
{2}=22...,$ then $\Pi(\mathbf{i})$ is also tiling by half unit intervals, but
the support is either $[k,\infty)$ or $(-\infty,k]$ for some $k\in\mathbb{Z}$.
\end{example}

\begin{example}
Let $F=\mathbb{\{R}^{2};f_{1},f_{2}\},f_{1}(x,y)=(\frac{x}{2},\frac{y}%
{2}),f_{2}(x,y)=(\frac{x+1}{2},\frac{y}{2}).$ The attractor is $A=[0,1]\times
\{0\},$ and the central open set is $C=(0,1)\times(-\infty,\infty)$ is
unbounded. If the tail of $\mathbf{i\in}\Sigma$ is neither $\overline
{1}=11...$ nor $\overline{2}=22...,$ then $A=[0,1]\times\{0\}$, $\Pi
(\mathbf{i})=$ $\{[\frac{n}{2},\frac{n+1}{2}]\times\{0\}\subset\mathbb{R}%
^{2}|n\in\mathbb{Z}\},$ and $\Xi(\mathbf{i})=$ $\{[\frac{n}{2},\frac{n+1}%
{2}]\times(-\infty,\infty):n\in\mathbb{Z}\}$.
\end{example}

\begin{figure}[ptb]
\centering
\includegraphics[
height=1.345in,
width=5.1814in
]{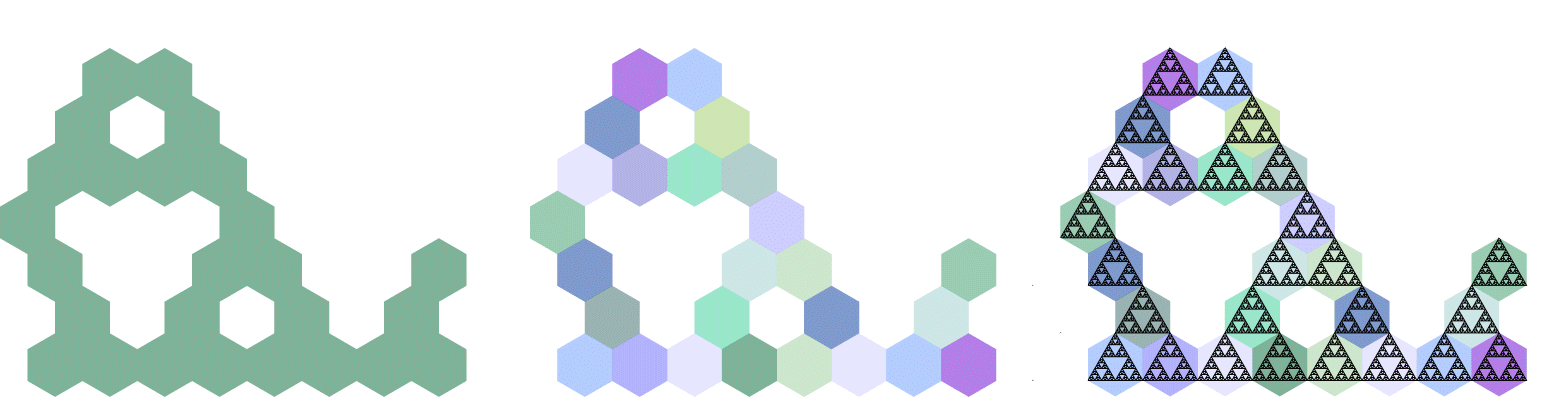}\label{sierptiles}\caption{From left to right: (i) part of
the support of a central open set tiling; (ii) a patch of the same central
open set tiling; (iii) same patch, also showing the underlying fractal tiles.
See Example \ref{sierptiles1ex}.}%
\label{sierptiles1}%
\end{figure}

\begin{example}
\label{sierptiles1ex} Let $F=\{\mathbb{R}^{2};f_{1,}f_{2},f_{3}\}$ where
\[
f_{1}=%
\begin{bmatrix}
\frac{1}{2} & 0 & 0\\
0 & \frac{1}{2} & 0
\end{bmatrix}
,f_{2}=%
\begin{bmatrix}
\frac{1}{2} & 0 & \frac{1}{2}\\
0 & \frac{1}{2} & 0
\end{bmatrix}
,f_{3}=%
\begin{bmatrix}
\frac{1}{2} & 0 & \frac{1}{4}\\
0 & \frac{1}{2} & \frac{\sqrt{3}}{4}%
\end{bmatrix}
.
\]
Then the attractor is a Sierpinski triange and the central open set is a
hexagon, see \cite{bandtneighbor2}. See Figure \ref{sierptiles1}.
\end{example}

\begin{figure}[ptb]
\centering
\includegraphics[
height=5.4744in,
width=4.8468in
]{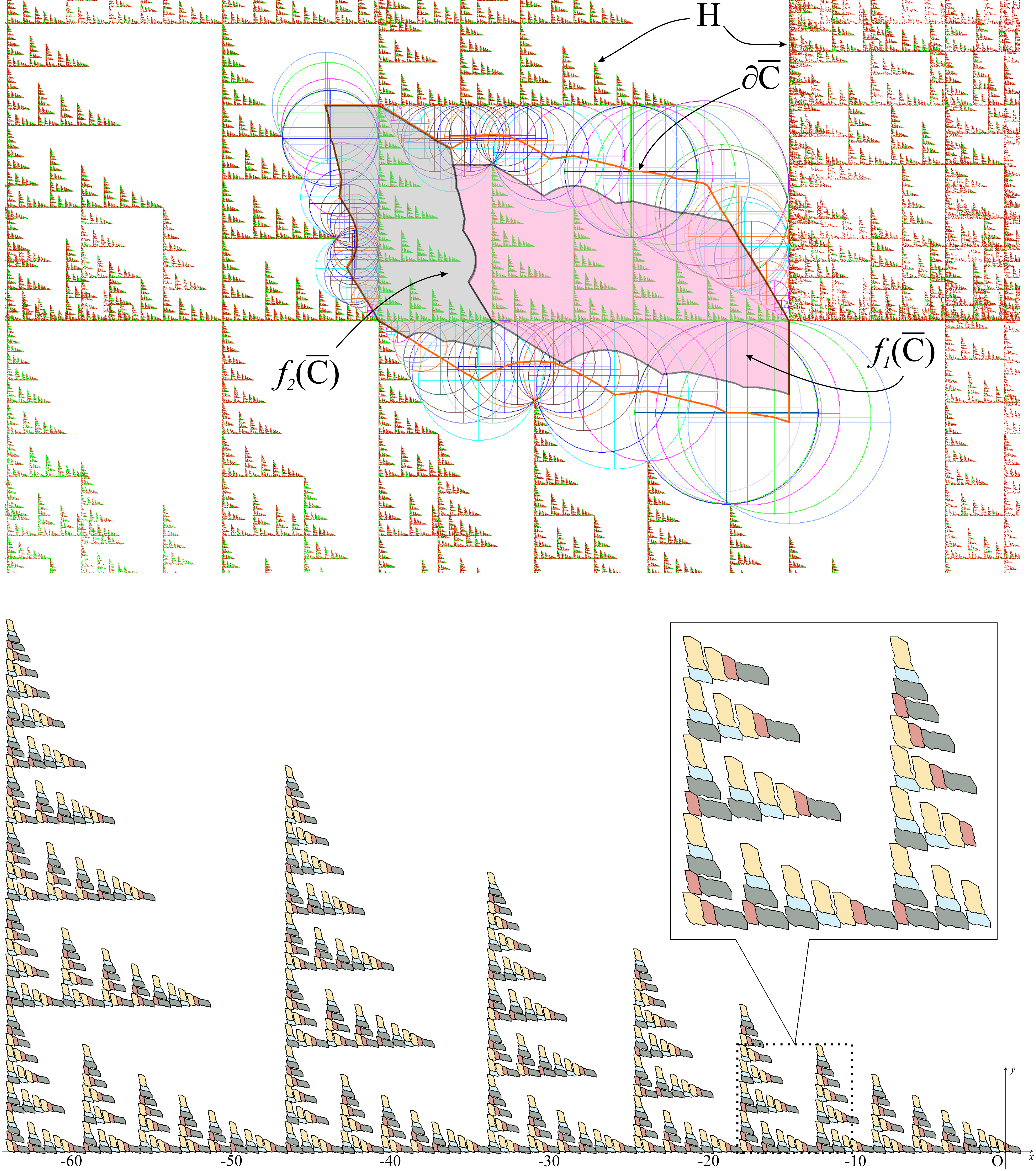}\caption{The top panel illustrates (i) $H$ for the IFS in
Example \ref{combplusex}; (ii) $\partial\overline{C},$ the boundary of the
central open set $C$; (iii) $f_{1}(\overline{C})$ and $f_{2}(\overline{C});$
(iv) circles whose centers approximate points on $\partial C$. The bottom
panel illustrates $\Xi(1111111111111),$ and (inset) a patch of this tiling
that also appears in Figure \ref{exactsbys}. }%
\label{combplus4}%
\end{figure}In Theorem \ref{maintheorem} below we establish some properties of
the collection of tilings $\mathbb{T}_{T}$. We use the following terminology.
We say that $\mathbf{i}\in\Sigma^{\infty}$ is \textit{disjunctiv}e, when given
any finite word $\mathbf{j}|p=j_{1}j_{2}\dots j_{p}$, there is $k\in
\mathbb{N}$ such that $i_{k+1}\dots i_{k+p}=j_{1}j_{2}\dots j_{p}.$ We say
that $\mathbf{i}\in\Sigma^{\infty}$ is \textit{reversible} when $A$ has
non-empty interior $A^{\circ}$ and there exists $k<l$ such that $f_{l}^{{}%
}f_{l-1}\dots f_{k}(A)\subset A^{\circ}$, the interior of $A$. If $A$ has
non-empty interior, then disjunctive is a special case of reversible,
\cite{barnvince2013, tilings}, but disjunctiveness is much easier to check
than reversibility.

\begin{theorem}
\label{maintheorem} Let $F=\{\mathbb{R}^{n};f_{1},\dots f_{m}\}$ be an IFS of
contractive similitudes. Let $T\subset\mathbb{R}^{n}$ be closed and let $c$ be
a cost function. Then $\Pi_{T}(\mathbf{i}|k)$ and $\Pi_{T}(\mathbf{i})$ are
well-defined collections of closed subsets of $\mathbb{R}^{n}$, (i.e. they are
tilings), and Equation (\ref{nestedeq}) holds for all $k\in\mathbb{N}$ and all
$\mathbf{i\in}\Sigma^{\infty}$.

If $F$ obeys the OSC, then the following statements are true.

(i) For all $\mathbf{i}\in\Sigma^{\infty}$ the interiors of the tiles that
comprise $\Xi(\mathbf{i})=\Pi_{\overline{C}}(\mathbf{i})$ are disjoint.

(ii) The interior of $f_{(\mathbf{i|}k\mathbf{)}}^{-1}f_{(\mathbf{j|}%
l\mathbf{)}}(\overline{C})$ is the central open set for the iterated function
system $SFS^{-1}$ where $S=f_{(\mathbf{i|}k\mathbf{)}}^{-1}f_{(\mathbf{j|}%
l\mathbf{)}}$ is a similitude. In this sense all tiles in $\Xi(\mathbf{i})$
are central open sets.

(iii) In $\mathbb{R}^{2}$, if $A$ is a polygon and $c_{i}=a_{i}\in\mathbb{N}$
for all $i\in\left\{  1,...m\right\}  $, then $\Xi(\mathbf{i})=\Pi
(\mathbf{i})$. In this case, if also $\mathbf{i}\in\Sigma^{\infty}$ is
reversible, then the support of $\Pi(\mathbf{i})$ is $\mathbb{R}^{n}$ and
$\Pi(\mathbf{i})$ is a self-similar tiling in the sense of standard works such
as \cite{anderson, grun, solomyak} and many others such as \cite{frank2,
kenyon, lagarias, radin} $.$

(iv) In general the tiling $\Pi_{T}(\mathbf{i})\ $is incommensurate both as
defined here and in the sense of \cite{smilansky}. The tilings $\Pi
_{T}(\mathbf{i}),$ in particular $\Pi(\mathbf{i})$ and $\Xi(\mathbf{i}),$ are
commensurate when $c_{i}=a_{i}\in\mathbb{N}$ for all $i=1,2,\dots m$.

(v) Let $\mathbf{i,j}$ $\in\Sigma^{\infty}$, $p,q\in\mathbb{N}$, and the cost
function $c,$ be such that $\sigma^{p}\mathbf{i}=\sigma^{q}\mathbf{j}$,
$c(\mathbf{i}|p)=c(\mathbf{j|}q)$ and $c_{l}=a_{l}$ for all $l=\{1,2,\dots
m\}$. Then
\begin{equation}
\Pi_{T}(\mathbf{i)=}E\Pi_{T}(\mathbf{j)}\label{rigidequation}%
\end{equation}
for all $T$, where $E=f_{i_{1}}^{-1}\dots f_{i_{p}}^{-1}f_{j_{q}}\dots
f_{j_{1}}$.

(vi) If $\mathbf{i}$ $\in\Sigma^{\infty}$ is reversible and $A$ has non-empty
interior, then the support of $\Pi(\mathbf{i),}$ namely $\bigcup\{t\in\Pi
_{A}(\mathbf{i)}\},$ is $\mathbb{R}^{n}$.
\end{theorem}

\begin{figure}[ptb]
\centering\includegraphics[
height=4.0299in,
width=5.2079in
]{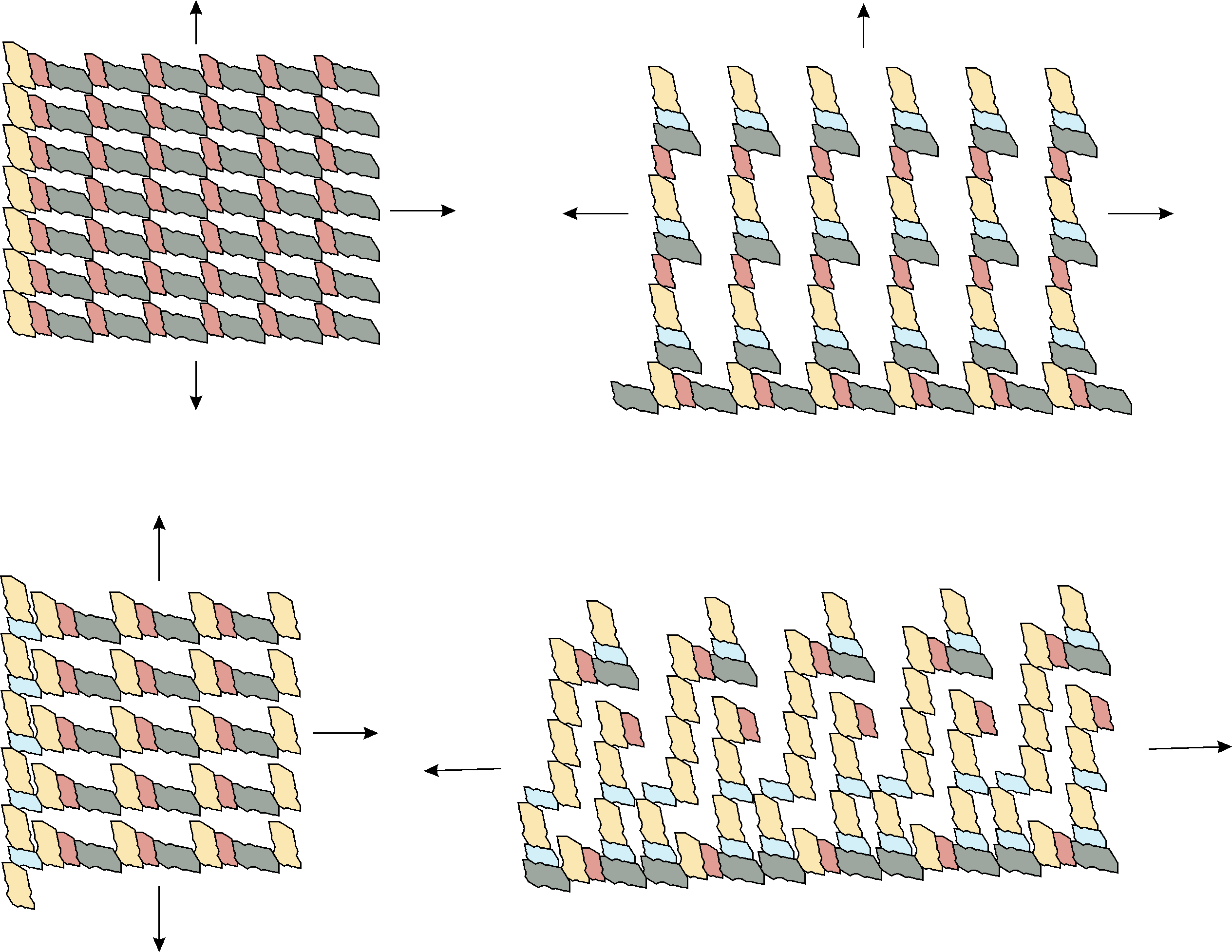}\caption{Examples of some tiling patterns that can be
formed using the prototiles in Figure \ref{prototiles4}, when the tiles are
allowed to touch only as in Figure \ref{prototiles4}. The arrows show
directions in which part of the tiling may be repeated periodically.}%
\label{combtilefill4}%
\end{figure}

\begin{figure}[ptb]
\centering
\includegraphics[
height=4.1635in,
width=4.5247in
]{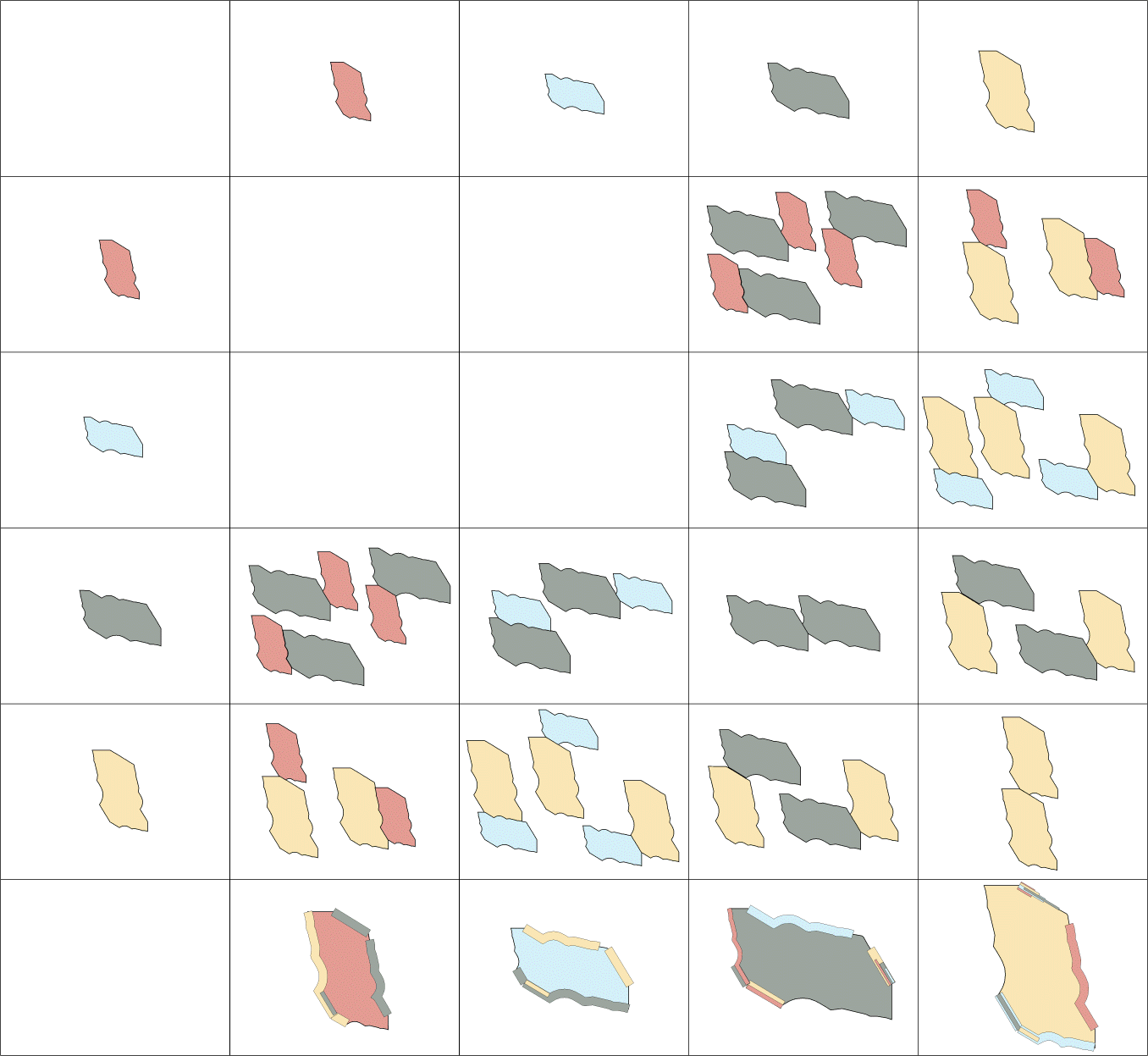}\caption{Left column: The four prototiles, with respect to
translations, in Figures \ref{exactsbys}, \ref{combplus4}, \ref{combtilefill4}%
. The bottom row illustrates the points of contact in every central open set
tiling $\Xi(\mathbf{i})$, $\mathbf{i}\in\Sigma^{\infty}$. The remaining boxes
illustrate all of the allowed combinations between pairs of tiles, indexed by
row and column.}%
\label{prototiles4}%
\end{figure}

\begin{proof}
The initial assertion follows at once from Equation (\ref{nestedeq}). This
generalizes, in the case of a single vertex, a core result in
\cite{barnvince2020}$.$

(i) We need to show that for all $\mathbf{i}\in\Sigma^{\infty}$ the interiors
of the tiles that comprise $\Xi(\mathbf{i})=\Pi_{\overline{C}}(\mathbf{i})$
are disjoint. Suppose
\[
f_{(\mathbf{i|}k\mathbf{)}}^{-1}f_{(\mathbf{j|}l\mathbf{)}}^{{}}(C)\cap
f_{(\mathbf{i|}p\mathbf{)}}^{-1}f_{(\mathbf{t|}q\mathbf{)}}(C)\neq\emptyset
\]
for some $k,p,l,q\in\mathbb{N}$, $\mathbf{j,t}\in\Sigma^{\infty}$ such that%
\[
c(\mathbf{j}|l-1)\leq c(\mathbf{i}|k)<c(\mathbf{j}|l)\text{ and }%
c(\mathbf{t}|q-1)\leq c(\mathbf{i}|p)<c(\mathbf{t}|q).
\]
We can assume $k<p$, $i_{k}\neq j_{1}$ and $i_{p}\neq t_{1}$. It follows that%
\[
f_{j_{1}}...f_{j_{l}}(C)\cap f_{i_{k+1}}^{-1}\dots f_{i_{p}}^{-1}f_{t_{1}%
}...f_{t_{q}}(C)\neq\emptyset.
\]
This implies%
\[
f_{i_{p}}^{{}}\dots f_{i_{k+1}}^{{}}f_{j_{1}}...f_{j_{l}}(C)\cap f_{t_{1}%
}...f_{t_{q}}(C)\neq\emptyset.
\]
This implies%
\[
f_{i_{p}}^{{}}(C)\cap f_{t_{1}}(C)\neq\emptyset,
\]
where $i_{p}\neq t_{1}$ which contradicts the fact that $C$ obeys the OSC.

(ii) This is an exercise in change of coordinates. We show that the interior
of $f_{(\mathbf{i|}k\mathbf{)}}^{-1}f_{(\mathbf{j|}l\mathbf{)}}(\overline{C})$
is the central open set for the iterated function system $SFS^{-1}$ where
$S=f_{(\mathbf{i|}k\mathbf{)}}^{-1}f_{(\mathbf{j|}l\mathbf{)}}$ is a similitude.

Let $S$ be the similitude $f_{(\mathbf{i|}k\mathbf{)}}^{-1}f_{(\mathbf{j|}%
l\mathbf{)}}.$ Then $A^{\prime}:=f_{(\mathbf{i|}k\mathbf{)}}^{-1}%
f_{(\mathbf{j|}l\mathbf{)}}(\overline{C})=SA$ is the attractor of the IFS
$F^{\prime}=S^{{}}FS^{-1}:=\{\mathbb{R}^{n};Sf_{i}S^{-1}|i=1,2,\dots m\}.$ The
neighbor maps of $F^{\prime}$ are $\left\{  Sf_{\mathbf{p}}^{-1}f_{\mathbf{q}%
}^{{}}S^{-1}|p_{1}\neq q_{1},\mathbf{p,q\in}\Sigma^{\ast}\right\}  $, so the
central open set for $F^{\prime}$ is
\begin{align*}
C^{\prime} &  =\{x^{\prime}\in\mathbb{R}^{m}|d(x^{\prime},A^{\prime
})<d(x,C^{\prime})\}\\
&  =\{x^{\prime}\in\mathbb{R}^{m}|d(x^{\prime},SA)<d(x^{\prime},SC)\}\\
&  =\{Sx\in\mathbb{R}^{m}|d(Sx,SA)<d(Sx,SC)\}\\
&  =S\{x\in\mathbb{R}^{m}|d(x,A)<d(x,C)\}=SC
\end{align*}
where in the penultimate step we have used the fact that $S$ is a similitude.

(iii) This follows from \cite{polygon}, which considers the case of
self-similar polygonal tilings, upon noting that if $A$ has non-empty interior
and obeys the OSC then $A=\overline{C}$ is a polygon.

(iv) If $c_{i}=a_{i}\in\mathbb{N}$ for all $i$, then it is readily seen that
each tile in $\Pi(\mathbf{i})$ is a copy of $T$ scaled by $s^{a}$ for some
$a\in\{1,2,\dots a_{\max}\}$ where $a_{\max}=\max\{a_{1},a_{2},\dots a_{m}\}.$

(v) Since $i_{p+1}i_{p+2}\dots=j_{q+1}j_{q+2}\dots$it follows that
\[
\Pi_{T}(\sigma^{p}\mathbf{i)=}\Pi_{T}(\sigma^{q}\mathbf{j)}%
\]
where $\sigma:\Sigma^{\infty}\rightarrow\Sigma^{\infty}$ is the shift
operator. Since $c(\mathbf{i}|p)=c(\mathbf{j|}q)$, we have
\begin{align*}
\{f_{j_{1}^{\prime}}...f_{j_{l}^{\prime}}(T)|\mathbf{j}^{\prime} &  \in
\Sigma^{\infty},l\in\mathbb{N},c(\mathbf{j}^{\prime}|l-1)\leq c(\mathbf{i}%
|p)<c(\mathbf{j}^{\prime}|l)\}\\
&  =\\
\{f_{j_{1}}...f_{j_{l}}(T)|\mathbf{j}^{\prime} &  \in\Sigma^{\infty}%
,l\in\mathbb{N},c(\mathbf{j}^{\prime}|l-1)\leq c(\mathbf{j|}q)<c(\mathbf{j}%
^{\prime}|l)\}.
\end{align*}
The result now follows from
\begin{align*}
f_{i_{p}}^{{}}\dots f_{i_{1}}^{{}}\bigcup\limits_{k\geq p}\{f_{i_{1}}%
^{-1}\dots f_{i_{k}}^{-1}(\{f_{j_{1}}...f_{j_{l}}(T)|\mathbf{j}^{\prime} &
\in\Sigma^{\infty},l\in\mathbb{N},c(\mathbf{j}|l-1)\leq c(\mathbf{i}%
|k)<c(\mathbf{j}|l)\})\}\\
&  =\\
f_{j_{q}}\dots f_{j_{1}}\bigcup\limits_{k\geq q}\{f_{j_{1}}^{-1}\dots
f_{j_{k}}^{-1}(\{f_{j_{1}}...f_{j_{l}}(T)|\mathbf{j}^{\prime} &  \in
\Sigma^{\infty},l\in\mathbb{N},c(\mathbf{j}|l-1)\leq c(\mathbf{i}%
|k)<c(\mathbf{j}|l)\})\}.
\end{align*}

(vi) This follows similar lines to \cite{tilings, barnvince2020} and is
omitted here. The argument there rests on the observation that, for reversible
addresses $\mathbf{i}$, the support of the tiling $\Pi(\mathbf{i|}k\mathbf{)}$
is contained in the interior of the tiling $\Pi(\mathbf{i|}k+l\mathbf{)}$ for
large enough $l$.
\end{proof}

When the IFS $F$ is rigid\textit{ }and $c_{i}=a_{i}\in\mathbb{N}$ for all
$i=1,2,\dots m,$ a converse of (v) in the Theorem is true. We say that $F$ is
\textit{rigid} (with respect to translations) when the statement
\textquotedblleft$T_{k}$ meets $ET_{l}$\textquotedblright\ for any
$k,l\in\{1,2,\dots\max\{a_{i}\}\}$ and any translation $E,$ implies
\textquotedblleft$T_{k}\ $is contained in $ET_{l}$ or
vice-versa\textquotedblright. Both examples in Figure \ref{leaftiles4} are
rigid. The only way that a rigid tiling $\Pi(\mathbf{i}$) can meet a
translation of another rigid tiling $\Pi(\mathbf{j})$ is when Equation
\ref{rigidequation} holds. Such tilings cannot be periodic and have
interesting properties, see \cite{barnvince2020} and references.

\section{Calculations and examples}

In this Section we present examples, including ones which show how we
calculate approximations to central open sets.

\begin{example}
\label{combplusex} We consider the IFS $\{\mathbb{R}^{2};f_{1},f_{2}\}$ where
\[
f_{1}=%
\begin{bmatrix}
s & 0 & 1\\
0 & s & 0
\end{bmatrix}
,f_{2}=%
\begin{bmatrix}
0 & s^{2} & 0\\
s^{2} & 0 & 0
\end{bmatrix}
\]
with $s+s^{4}=1$, $s>0$. This IFS was mentioned in Example \ref{fig2-introex}
and its attractor is illustrated at the bottom left in Figure \ref{fig2-intro}
. In the top image in Figure \ref{combplus4} we illustrate the central open
set and suggest how it was approximated. In this example and others the
attractor $A$ and parts of the fast basin $B$ were calculated by using random
iteration \cite{chaos}. The portion of $B\backslash A$ closest to $A$ was
assumed to be the union of the sets $f_{\mathbf{i}}^{-1}f_{\mathbf{j}}^{{}%
}(A)$ for $i_{1}\neq j_{1}$ and $\mathbf{i,j}\in\Sigma^{4},$ also computed by
random iteration. In order to estimate points on $\partial C=\partial
\overline{C},$ circles that appeared to touch both $B\backslash A$ and $A$
were constructed$.$ Calculations and constructions were performed on digital
images of resolution $1024\times1024$. The bottom image in Figure
\ref{combplus4} illustrates $\Xi(1111111111111).$ This picture was constructed
by starting from a computed image of $\Pi(1111111111111)$ which is a
translation of the canonical tiling $T_{13}.$ Figure \ref{combtilefill4}
illustrates four tiling patterns constructed using the tiling rules in Figure
\ref{prototiles4}. The arrows point in directions in which a part of the
pattern could be repeated periodically.
\end{example}

\begin{figure}[ptb]
\centering
\includegraphics[
height=5.5865in,
width=4.9423in
]{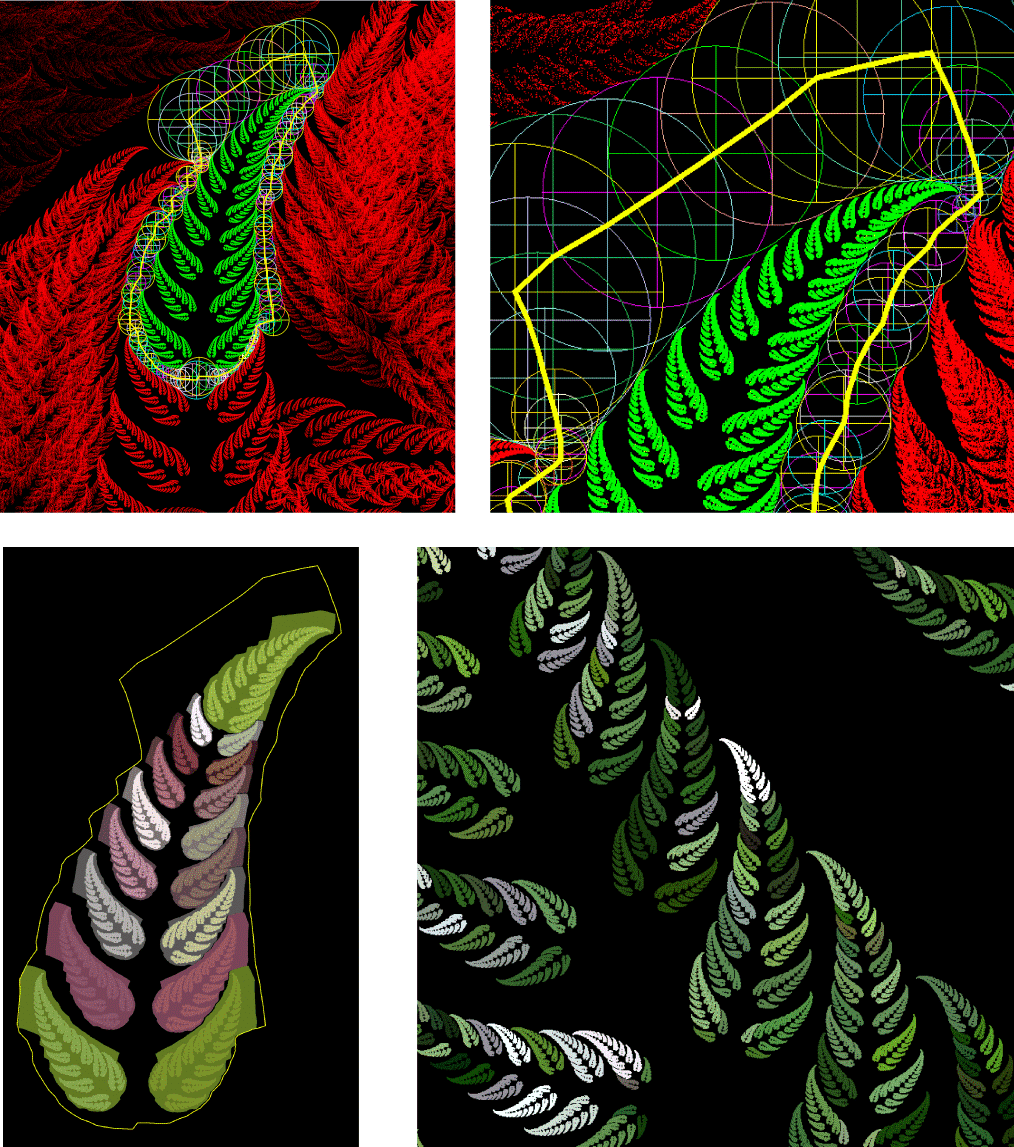}\caption{The top panels illustrate part of $B/A$ where $B$ is
the fast basin and $A$ is the attractor, for the IFS in Example
\ref{ferncomboex}; $\partial C,$ the boundary (yellow) of the central open set
$C$;  $\ $and circles whose centers approximate points on $\partial C$. The
bottom left panel shows scaled copies of $C$ in part of a tiling of the form
$\Xi(11111111)$ overlayed on $\Pi(11111111)$. The bottom right panel shows
part of  $\Pi(\mathbf{i})$ for some $\mathbf{i\in}\Sigma^{\infty}$.}%
\label{ferncombo}%
\end{figure}

\begin{example}
\label{ferncomboex} We consider the IFS $F=\{\mathbb{R}^{2};f_{1},f_{2}%
,f_{3}\}$ defined by%
\begin{equation}
f_{1}=%
\begin{bmatrix}
.85 & -.05 & .53842\\
.05 & .85 & -.15789
\end{bmatrix}
,f_{2}=%
\begin{bmatrix}
.17 & .22 & .195909\\
-.22 & .17 & .776364
\end{bmatrix}
,f_{3}=%
\begin{bmatrix}
-.17 & -.22 & .805\\
-.22 & .17 & .776364
\end{bmatrix}
.\label{fernifs}%
\end{equation}
The attractor is illustrated in green in the central zone of the top left
panel in Figure \ref{ferncombo}. The attractor is totally disconnected
although the image makes it appear to have connected components, because of
digitization effects. The red set in the top left panel is an approximation to
the relevant part of $B\backslash A$ and was calculated in the same way as in
Example \ref{combplusex}. A close-up is shown on the right, illustrating how
we estimated the central open set. To make pictures of some tilings we chose
$c_{1}=1$ and $c_{2}=c_{3}=8.$ The scaling factor for $f_{1}$ is
$s_{1}=\allowbreak0.851\dots$ while the scaling factors for $f_{2}$ and
$f_{3}$ are $s_{2}=s_{3}=0.278\dots$ so $s_{1}^{8}\simeq s_{2}=s_{3},$ but
$s_{1}^{8}\neq s_{2}$. Thus $\Pi(\mathbf{i)}$ and $\Xi(\mathbf{i)}$
incommensurate tilings for any $\mathbf{i\in}\Sigma^{\infty}$. The bottom left
image illustrates the tiling $\Xi(1111111)$ surrounded by $f_{1}^{-8}\partial
C$. It illustrates the relationship between $\partial C$ and some tiles. We
observed that there appeared to be eight different tile sizes in any square
patch of $\Pi(\mathbf{i})$ digitized at resolution $2048\times2048,$ when
keeping $s_{1}^{-8}A$ to be roughly the size of the viewing window. This
accords with a comment in \cite{bandtneighbor2} regarding a result of Schief
\cite{schief}: \textquotedblleft\textit{There exists an integer }$N$\textit{
such that at most }$N$\textit{ incomparable pieces }$A_{\mathbf{j}}$%
\textit{(}$=f_{\mathbf{j}}(A)$\textit{) of size }$\geq\varepsilon$\textit{ can
intersect the }$\varepsilon$\textit{-neighborhood (sic) of a piece
}$A_{\mathbf{i}}$\textit{ of diameter }$\varepsilon$\textit{.}%
\textquotedblright\ (The sets $A_{j_{1}\dots j_{n}}$ and $A_{i_{1}\dots i_{m}%
}$, referred to in the quote as \textquotedblleft pieces\textquotedblright,
are said to be incomparable if there exists no $k_{1}\dots k_{p}$ such that
$j_{1}\dots j_{n}=i_{1}\dots i_{m}k_{1}\dots k_{p}$ or $i_{1}\dots i_{m}%
=j_{1}\dots j_{n}k_{1}\dots k_{p}$.) An example of part an unbounded tiling
$\Pi(\mathbf{i})$ is illustrated at the bottom right in Figure \ref{ferncombo}%
$.$
\end{example}

\begin{figure}[ptb]
\centering
\includegraphics[
height=4.2225in,
width=5.0104in
]{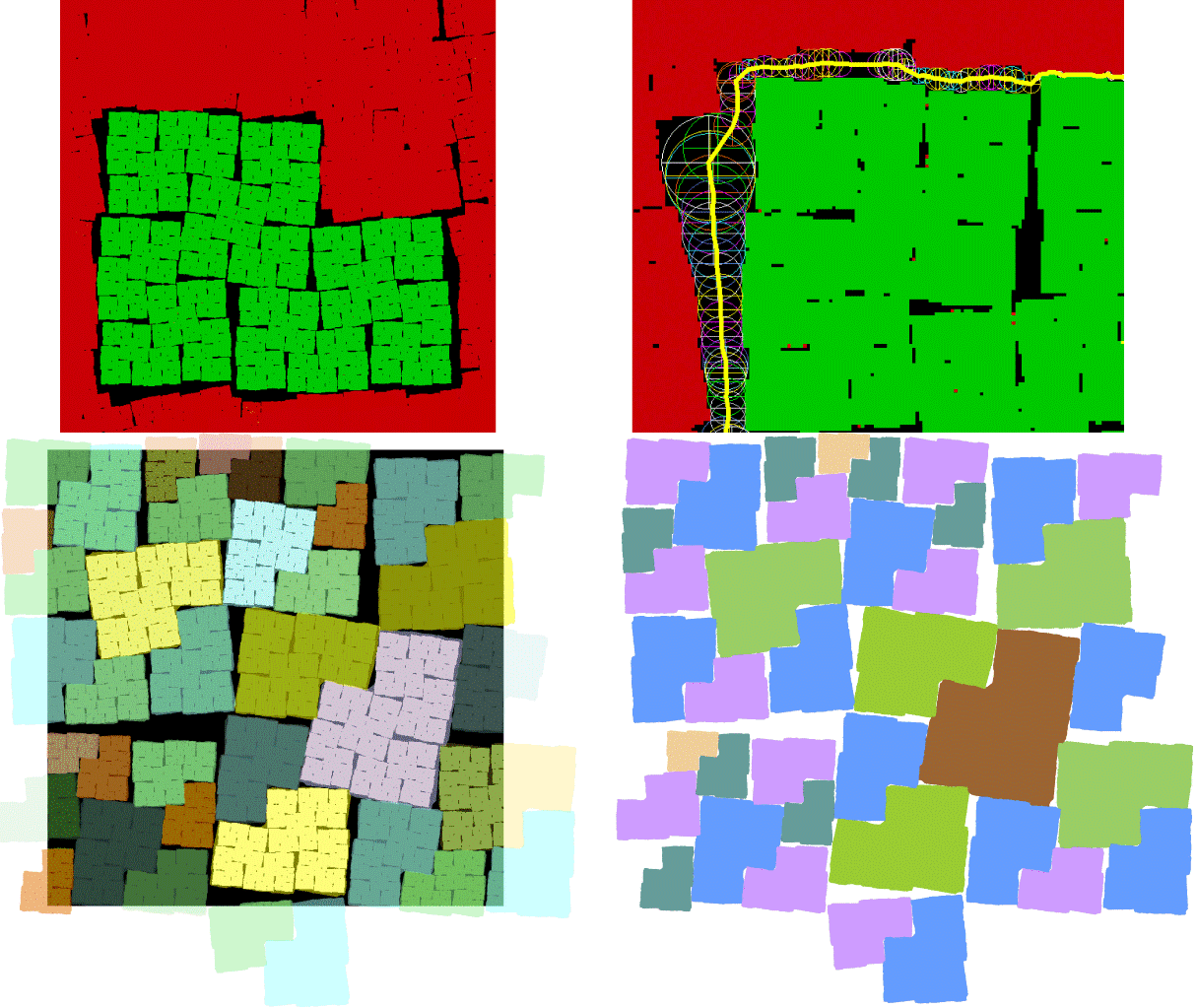}\caption{See Example \ref{combinedex}. Top left illustrates the
attractor of an IFS of two maps (green), and its fast basin minus the
attractor (red). The top right hand image illustrates how the boundary, shown
in yellow, of the central open set $C$ was estimated. The lower left image
shows a patch of $\Pi(\mathbf{i})$, specifically the tiles which meet a (dark)
square. The lower right image shows the corresponding patch of  of
$\Xi(\mathbf{i}).$}%
\label{combined}%
\end{figure}

\begin{example}
\label{combinedex} In Figure \ref{combined} we illustrate the calculation of
the central open set for an IFS involving two maps. The method is the same as
described in connection with Figures \ref{ferncombo} and \ref{combplus4}. Here
the IFS is close to the one illustrated in the top row of Figure
\ref{fig2-intro}, see Example \ref{fig2-introex}. Each map here is slightly
more contractive and rotated by a small amount. The IFS here is $\{\mathbb{R}%
^{2};f_{1},f_{2}\}$ where
\[
f_{1}=%
\begin{bmatrix}
-.02447 & .777910 & 0\\
-.77791 & -.02447 & .78615
\end{bmatrix}
,f_{2}=%
\begin{bmatrix}
.61156 & -.019221 & -.019221\\
-.019221 & .61156 & 0
\end{bmatrix}
\]
As in Example \ref{ferncomboex} this numerical model is only approximately
scaling, but we treat it as though it is by choosing $c_{1}=1$ and $c_{2}=2.$
We find this example interesting because it suggests natural situations
involving cracks.
\end{example}

\begin{figure}[ptb]
\centering
\includegraphics[
height=2.2383in,
width=4.8858in
]{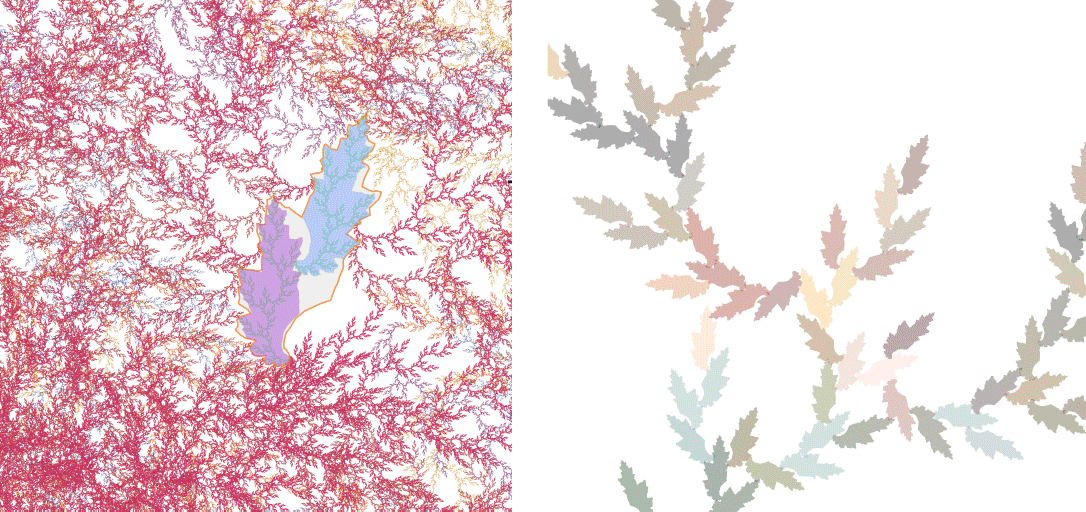}\caption{See Example \ref{newgrowthsbysex}. The left image
shows an approximate central open set $\protect\widetilde{C}$ and
corresponding fast basin. The sets $f_{1}(\protect\widetilde{C})$ and
$f_{2}(\protect\widetilde{C})$ are shown to fit neatly inside
$\protect\widetilde{C}$, and are surrounded by an approximation to
$B\backslash A$. It may well be true that the attractor is overlapping and
there is no central open set. Nonetheless, a tiling $\Pi
_{\protect\widetilde{C}}(\mathbf{i}),$ a patch of which is illustrated on the
right, is only slightly overlapping.}%
\label{newgrowthsbys}%
\end{figure}

\begin{example}
\label{newgrowthsbysex} See Figure \ref{newgrowthsbys}. The IFS here is
$\{\mathbb{R}^{2};f_{1},f_{2}\}$ where
\[
f_{1}=%
\begin{bmatrix}
.6413 & -.3283 & .3231\\
.3283 & .6413 & -.133
\end{bmatrix}
,f_{2}=%
\begin{bmatrix}
-.2362 & .4620 & .8052\\
.4620 & .2362 & .5093
\end{bmatrix}
\]
One can see how $\overline{C}$ is approximately tiled by $f_{1}(\overline{C})$
and $f_{2}(\overline{C})$ in the left image. That is, one can see the
relationship between $\overline{C}$ and $\left\{  f_{1}(\overline{C}
),f_{2}(\overline{C})\right\}  $. Note that the actual attractor may not obey
the OSC, our estimated \textquotedblleft central open set\textquotedblright%
\ $\widetilde{C}$ may not obey the OSC, and the tiling $\Pi_{\widetilde{C}%
}(\mathbf{i})$ may be overlapping.
\end{example}

We thank to Krystof Lesniak and Nina Snigreva for many helpful comments and suggestions.

\vspace{0.5cm}

\noindent\makebox[\linewidth]{\rule{80mm}{0.4pt}}\vspace{0.5cm}

\begin{center}
\noindent\begin{minipage}{8cm}
\centering
\small{
Mathematical Sciences Institute, Australian National University, Canberra, Australia\\
{\tt LouisaBarnsley@gmail.com}\vspace{0.2cm}
\noindent \\
{\tt Michael.Barnsley@anu.edu.au}}
\end{minipage}

\end{center}

\end{document}